\documentclass{article}
\usepackage {amsthm, amsmath, amssymb}
\newtheorem{theorem}{Theorem}
\newtheorem{corollary}[theorem]{Corollary}
\theoremstyle{remark}
\newtheorem{remark}[theorem]{Remark}

\newcommand {\calk} {\mathcal {K}}
\newcommand {\BS}[1]{\mathop\mathrm{BS}({#1})}
\newcommand {\BSL} {\mathrm{BS}_{\mathrm{lower}}}
\newcommand {\BSU} {\mathrm{BS}_{\mathrm{upper}}}
\newcommand {\Gal}{\mathop\mathrm{Gal}}

\title{The Brauer-Siegel and Tsfasman-Vl\u{a}du\c{t} Theorems for Almost Normal Extensions of Number Fields}
\author {Alexei Zykin\thanks{Independent University of Moscow, 119002, Bolshoy Vlasyevskiy 11, Moscow, Russia. E-mail:
zykin@mccme.ru}}
\date {}

\begin{document}
\maketitle
\begin{abstract}
The classical Brauer-Siegel theorem states that if $k$ runs through the sequence of normal extensions of $\mathbb{Q}$ such that $n_k/\log|D_k|\to 0,$ then $\log h_k R_k/\log \sqrt{|D_k|}\to 1.$ First, in this paper we obtain the generalization of the Brauer-Siegel and Tsfasman-Vl\u{a}du\c{t} theorems to the case of almost normal number fields. Second, using the approach of Hajir and Maire, we construct several new examples concerning the Brauer-Siegel ratio in asymptotically good towers of number fields. These examples give smaller values of the Brauer-Siegel ratio than those given by Tsfasman and Vl\u{a}du\c{t}. 
\end{abstract}
\section{Introduction}
Let $K$ be an algebraic number field of degree $n_K=[K: \mathbb{Q}]$ and discriminant $D_K.$ We define the genus of $K$ as $g_K=\log\sqrt{D_K}.$ By $h_K$ we denote the class-number of $K,$ $R_K$ denotes its regulator. We call a sequence $\{K_i\}$ of number fields a family if $K_i$ is non-isomorphic to $K_j$ for $i\ne j.$ A family is called a tower if also $K_i\subset K_{i+1}$ for any $i.$ For a family of number fields we consider the limit
$$\BS{\calk}:=\lim_{i\to\infty}\frac{\log h_{K_i}R_{K_i}}{g_{K_i}}.$$
The classical Brauer-Siegel theorem, proved by Brauer (see~\cite{Bra}), states that for a family $\calk=\{K_i\}$ we have $\BS{\calk}= 1$ if the family satisfies two conditions:
\begin{itemize}
\item[(i)] $\lim\limits_{i\to \infty}\frac {n_{K_i}}{g_{K_i}}=0;$
\item[(ii)] either the generalized Riemann hypothesis (GRH) holds, or all the fields $K_i$ are normal over $\mathbb{Q}.$
\end{itemize}
 
We call a number field almost normal if there exists a finite tower of number fields $\mathbb{Q}=K_0\subset K_1\subset\dots\subset K_m=K$ such that all the extensions $K_i/K_{i-1}$ are normal. Weakening the condition (ii), we prove the following generalization of the classical Brauer-Siegel theorem to the case of almost normal number fields:
\begin{theorem}
\label{main}
Let $\calk=\{K_i\}$ be a family of almost normal number fields for which $n_{K_i}/g_{K_i}\to 0$ as $i\to\infty.$ Then we have $\BS{\calk}=1.$
\end{theorem}

It was shown by M. A. Tsfasman and S. G. Vl\u{a}du\c{t} that, taking in account non-archimedian places, one may generalize the Brauer-Siegel theorem to the case of extensions where the condition (i) does not hold. 

For a prime power $q$ we set $$N_q(K_i):=|\{v\in P(K_i): \mathop\mathrm{Norm}(v)=q\}|,$$ where $P(K_i)$ is the set of non-archimedian places of $K_i.$ We also put $N_{\mathbb{R}}(K_i)=r_1(K_i)$ and $N_{\mathbb{C}}(K_i)=r_2(K_i),$ where $r_1$ and $r_2$ stand for the number of real and (pairs of) complex embeddings. 

We consider the set $A = \{\mathbb{R}, \mathbb{C}; 2, 3, 4, 5, 7, 8, 9,\hdots\}$ of all prime powers plus two auxiliary symbols $\mathbb{R}$ and $\mathbb{C}$ as the set of indices. A family $\calk=\{K_i\}$ is called asymptotically exact if and only if for any $\alpha\in A$ the following limit exists:
$$\phi_{\alpha} = \phi_{\alpha}(\calk):=\lim_{i\to\infty}\frac{N_\alpha (K_i)}{g_{K_i}}.$$
We call an asymptotically exact family $\calk$ asymptotically good (respectively, bad) if there exists $\alpha\in A$ with $\phi_{\alpha}>0$ (respectively, $\phi_{\alpha}=0$ for any $\alpha\in A$). The condition on a family to be asymptotically bad is, in the number field case, obviously equivalent to the condition (i) in the classical Brauer-Siegel theorem. For an asymptotically good tower of number fields the following generalization of the Brauer-Siegel theorem was proved in~\cite{Tsfa}:
\begin{theorem}[Tsfasman-Vl\u{a}du\c{t} Theorem, see~\cite{Tsfa}, Theorem 7.3]
\label{assgood}
Assume that for an asymptotically good tower $\calk$ fields any of the following conditions is satisfied:
\begin{itemize}
\item GRH holds
\item All the fields $K_i$ are almost normal over $\mathbb{Q}.$
\end{itemize}
Then the limit $\BS{\calk} = \lim\limits_{i\to\infty}\frac{\log h_{K_i}R_{K_i}}{g_{K_i}}$ exists and we have:
\begin{equation}
\BS{\calk}=1+\sum_{q}\phi_q \log \frac{q}{q-1} - \phi_\mathbb{R} \log 2 - \phi_\mathbb{C}\log 2\pi,
\label{eqbs}
\end{equation}
the sum beeing taken over all prime powers $q$.
\end{theorem}
For an asymptotically bad tower of number fields we have $\phi_\mathbb{R} = 0$ and $\phi_\mathbb{C}=0$ as well as $\phi_q=0$ for all prime powers $q,$ so the right hand side of the formula~(\ref{eqbs}) equals to one. We also notice that the condition on a family to be asymptotically bad is equivalent to $\lim\limits_{i\to \infty}\frac {n_{K_i}}{g_{K_i}}=0.$ So, combining our theorem~\ref{main} with the theorem~\ref{assgood} we get the following corollary:
\begin{corollary}
For any tower $\calk=\{K_i\}$, $K_1\subset K_2\subset\hdots$ of almost normal number fields the limit $\BS{\calk}$ exists and we have:
$$\BS{\calk}=\lim_{i\to\infty}\frac{\log(h_i R_i)}{g_i}=1+\sum_{q}\phi_q \log \frac{q}{q-1} - \phi_\mathbb{R} \log 2 - \phi_\mathbb{C}\log 2\pi,$$
the sum beeing taken over all prime powers $q.$
\end{corollary}
In~\cite{Tsfa} bounds on the ratio $\BS{\calk}$ were given, together with examples showing that the value of $\BS{\calk}$ may be different from 1. We corrected some of these erroneous bounds and managed to precise a few of the estimates in the examples. Also, using the infinite tamely ramified towers, found by Hajir and Maire (see~\cite{Haj1}), we get (under GRH) new examples, both totaly complex and totally real, with the values of $\BS{\calk}$ smaller than those of totally real and totally complex examples of~\cite{Tsfa}. The result is as follows:
\begin{theorem}
\label{main1}
\begin{itemize}
\item[\em 1.]Let $k=\mathbb{Q}(\xi)$, where $\xi$ is a root of $f(x)=x^6+x^4-4x^3-7x^2-x+1,$ $K=k(\sqrt{\xi^5-467\xi^4+994\xi^3-3360\xi^2-2314\xi+961}).$ Then $K$ is totally complex and has an infinite tamely ramified $2$-tower $\calk$, for which, under GRH, we have: $$\BSL \leq \BS{\calk} \leq \BSU,$$
where $\BSL\approx0.56498\hdots$, $\BSU\approx0.59748\hdots$.
\item[\em 2.]Let $k=\mathbb{Q}(\xi)$, where $\xi$ is a root of $f(x)=x^6-x^5-10x^4+4x^3+29x^2+3x-13,$ $K=k(\sqrt{-2993\xi^5+7230\xi^4+18937\xi^3-38788\xi^2-32096\xi+44590}).$ Then $K$ is totally real and has an infinite tamely ramified $2$-tower $\calk$, for which, under GRH, we have: $$\BSL \leq \BS{\calk} \leq \BSU,$$
where $\BSL\approx0.79144\hdots$, $\BSU\approx0.81209\hdots$.
\end{itemize}
\end{theorem}
However, unconditionally (without GRH), the estimates for totally complex fields that may be obtained using the methods developed by Tsfasman and Vl\u{a}du\c{t} lead to slightly worse results, than those already known from~\cite{Tsfa}. This is due to a rather large number of prime ideals of small norm in the field $K.$ For the same reasons the upper bounds for the Brauer-Siegel ratio for other fields constructed in~\cite{Haj1} are too high, though the lower bounds are still good enough.

Finally we present the table (the ameliorated version of the table of~\cite{Tsfa}), where all the bounds and estimates are given together:
\footnotesize
\begin{center}
\begin{tabular}{|c|l|cccc|}
\hline
\multicolumn{2}{|c|}{}&lower&lower&upper&upper\\
\multicolumn{2}{|c|}{}&bound&example&example&bound\\
\hline
&all fields&0.5165&0.5649-0.5975&1.0602-1.0798&1.0938\\
GRH&totally real&0.7419&0.7914-0.8121&1.0602-1.0798&1.0938\\
&totally complex&0.5165&0.5649-0.5975&1.0482-1.0653&1.0764\\
\hline
&all fields&0.4087&0.5939-0.6208&1.0602-1.1133&1.1588\\
Unconditional&totally real&0.6625&0.8009-0.9081&1.0602-1.1133&1.1588\\
&totally complex&0.4087&0.5939-0.6208&1.0482-1.1026&1.1310\\
\hline
\end{tabular}
\end{center}
\normalsize
\section{Proof of Theorem~\ref{main}}
Let $\zeta_K(s)$ be the Dedekind zeta function of the number field $K$ and $\varkappa_K$ its residue at $s=1$. By $w_K$ we denote the number of roots of unity in $K,$ and by $r_1, r_2$ the number of real and complex places of $K$ respectively. We have the following residue formula (see ~\cite{Lan}, Chapter VIII, Section 3):
$$\varkappa=\frac{2^{r_1}(2\pi)^{r_2}h_KR_K}{w_K\sqrt{D_K}}.$$
Since
$$\sqrt{w_K/2}\leq\varphi(w_K)=[\mathbb{Q}\left(\zeta_{w_K}\right): \mathbb{Q}]\leq[K: \mathbb{Q}]=n_K,$$ we note that $w_K\leq 2n_K^2$ so $\log w_{K_j}/g_{K_j}\to 0.$ Thus, it is enough to prove that $\log\varkappa_{K_j}/\log D_{K_j}\to0.$

As for the upper bound we have
\begin{theorem}[See~\cite{Lou00}, Theorem 1]
\label{upperbound}
Let K be a number field of degree $n\geq2.$ Then, 
\begin{equation}
\label{estimlou1}
\varkappa_K\le\left(\frac{e\log D_K}{2(n-1)}\right)^{n-1}.
\end{equation}
Moreover, $1/2\leq\rho<1$ and $\zeta_K(\rho)=0$ imply
\begin{equation}
\label{estimlou2}
\varkappa_K\leq (1-\rho)\left(\frac{e\log D_K}{2n}\right)^n.
\end{equation}
\end{theorem}
Using the estimate~(\ref{estimlou1}) we get (even without the assumption of almost normality) the "easy inequality": 
$$\frac{\log\varkappa_{K_j}}{\log D_{K_j}}\le\frac{n_j-1}{\log D_{K_j}}\left(\log\frac e 2 + \log\frac{\log D_{K_j}}{n_j-1}\right)\to 0.$$

As for the lower bound the business is much more tricky and we will proceed to the proof after giving a few preliminary statements.

Let $K$ be a number field other than $\mathbb{Q}$. A real number $\rho$ is called an {\itshape exceptional zero} of $\zeta_K(s)$ if $\zeta_K(\rho)=0$ and $$1-(4\log D_K)^{-1}\leq\rho<1;$$
an exceptional zero $\rho$ of $\zeta_K(s)$ is called its {\itshape Siegel zero} if  
$$1-(16\log D_K)^{-1}\leq\rho<1.$$

Our proof will be based on the following fundamental property of Siegel zeroes proved by Stark:
\begin{theorem}[See~\cite{Sta}, Lemma 10]
\label{Zeroes}
Let $K$ be an almost normal number field, and let $\rho$ be a Siegel zero of $\zeta_K(s).$ Then there exists a quadratic subfield $k$ of $K$ such that $\zeta_k(\rho)=0.$
\end{theorem}

The next estimate is also due to Stark:
\begin{theorem}[See~\cite{Sta}, Lemma 4 or~\cite{Lou01}, Theorem 1]
\label{staineq}
Let $K$ be a number field and let $\rho$ be the exceptional zero of $\zeta_K(s)$ if it exists and $\rho=1-(4\log D_K)^{-1}$ otherwise. Then there is an absolute constant $c<1$ (effectively computable) such that 
\begin{equation}
\label{estimsta}
\varkappa_K>c(1-\rho)
\end{equation}
\end{theorem}

Our proof of Theorem~\ref{main} will be similar to the proof of the classical Brauer-Siegel theorem given in~\cite{Lou02}. We will use the Brauer-Siegel result for quadratic fields, a simple proof of which is given in~\cite{Gold}.
There are two cases to consider.
\begin{enumerate}
\item[1.] First, assume that $\zeta_{K_j}(s)$ has no Siegel zero. From~(\ref{estimsta}) we deduce that 
\begin{equation}
\label{firstcase}
\varkappa_{K_j}>c(1-\rho)\geq c\left(1-\left(1-\frac{1}{16\log D_{K_j}}\right)\right)=\frac{c}{16\log D_{K_j}}.
\end{equation}

\item[2.] Second, assume that there exists a Siegel zero $\rho$ of $\zeta_{K_j}(s).$ From Theorem~\ref{Zeroes} we see that there exists a quadratic subfield $k_j$ of $K_j$ such that  $\zeta_{k_j}(\rho)=0.$ Applying~(\ref{estimlou2}) and~(\ref{estimsta}) we obtain:
\begin{equation}
\label{secondcase}
\varkappa_{K_j}=\frac{\varkappa_{K_j}}{\varkappa_{k_j}}\varkappa_{k_j}\geq\frac{c(1-\rho)}{(1-\rho)\left(\frac{e\log D_{k_j}}{4}\right)^2}\varkappa_{k_j}=\frac{16c}{e^2\log^2D_{k_j}}\varkappa_{k_j}.
\end{equation}
\end{enumerate}

If the number of fields $K_j$ for which the second case holds is finite, then, using the fact that $\log D_{K_j}\to\infty,$ we get the desired lower estimate from~(\ref{firstcase}). 

Otherwise, we note that for a number field there exists at most one exceptional zero (See~\cite{Sta}, Lemma 3), so, applying this statement to the fields $k_j,$ we get that only finitely many of them may be isomorphic to each other and so $D_{k_j}\to\infty$ as $j\to\infty.$ Thus we may use the Brauer-Siegel result for quadratic fields: $$\frac{\log\varkappa_{k_j}}{\log D_{K_j}}\le\frac{\log\varkappa_{k_j}}{\log D_{k_j}}\to0.$$ Finally from~(\ref{secondcase}), we get:
$$\frac{\log\varkappa_{K_j}}{\log D_{K_j}}\geq \frac{16c}{e^2\log D_{K_j}}-2\frac{\log\log D_{k_j}}{\log D_{K_j}}+\frac{\log\varkappa_{k_j}}{\log D_{K_j}}\to0.$$
This concludes the proof. \qed
\begin{remark}
Our proof of Theorem~\ref{main} is explicit and effective if all the fields in the family $\calk$ contain no quadratic subfield and thus the corresponding zeta function does not have Siegel zeroes.
\end{remark}
\section{Proof of Theorem~\ref{main1}}
First we recall briefly some constructions related to class field towers. Let us fix a prime number $\ell.$ For a finitely generated pro-$\ell$ group $G$, we let $d(G)=\dim_{\mathbb{F}_\ell}H^1(G,\mathbb{F}_\ell)$ be its generator rank. Let $T$ be a finite set of ideals of a number field $K$ such that no prime in $T$ is a divisor of $\ell.$ We denote by $K_T$ the maximal $\ell$-extension of $K$ unramified outside $T, G_T=\Gal(K_T/K).$ We let
$$\theta_{K,T}=\begin{cases}
1, &\text{if $T\ne\varnothing$ and $K$ contains a primitive $\ell$th root of unity;}\\
0, &\text{otherwise.}
\end{cases}$$
Then we have (see~\cite{Shaf}, theorems 1 and 5):
\begin{theorem}
\label{inf}
If $d(G_T)\geq 2+2\sqrt{r_1(K)+r_2(K)+\theta_{K,T}},$ then $K_T$ is infinite.
\end{theorem}

To estimate $d(G_T)$ we use the following theorem 
\begin{theorem}[See~\cite{Mar}, section 2]
\label{GT}
Let $K/k$ be a finite Galois extension, $r_1=r_1(k)$, $r_2=r_2(k)$, $\rho$ be the number of real places of $k,$ ramified in $K, t$ be the number of primes in $k,$ ramified in $K$. We set $\delta_\ell = 1,$ if $k$ contains a primitive root of degree $\ell$ of unity and $\delta_\ell = 0$ otherwise. Then we have:
$$d(G_T)\geq d(G_\varnothing)\geq t-r_1-r_2+\rho-\delta_\ell$$
\end{theorem}

The number field arithmetic behind the construction of our theorem~\ref{main1} was mainly carried out with the help of the computer package PARI. However, we would like to present our examples in the way suitable for non-computer check. We give here the proof of the first part of our theorem, as the proof of the second part is very much similar and may be carried out simply by repeating all the steps of the proof given here.

We let $k=\mathbb{Q}(\xi),$ where $\xi$ is a root of $f(x)=x^6+x^4-4x^3-7x^2-x+1.$ Then $k$ is a field of signature $(4, 1)$ and discriminant $d_f=d_k=-23\cdot35509.$ Its ring of integers is $\mathcal{O}_k=\mathbb{Z}[\xi]$ and its class number is equal to 1. The principle ideal of norm $7\cdot13\cdot19^2\cdot23^2\cdot29\cdot31$ generated by $\eta=671\xi^5-467\xi^4+994\xi^3-3360\xi^2-2314\xi+961$ factors into eight different prime ideals of $\mathcal{O}_k.$ In fact, one may see that $\eta=\pi_7\pi_{13}\pi_{19}\pi'_{19}\pi_{23}\pi'_{23}\pi_{29}\pi_{31},$ where
\begin{align*}
&\pi_7=-9\xi^5+6\xi^4-13\xi^3+44\xi^2+31\xi-12,\\
&\pi_{13}=-7\xi^5+5\xi^4-11\xi^3+36\xi^2+23\xi-9,\\
&\pi_{19}=5\xi^5-4\xi^4+8\xi^3-26\xi^2-15\xi+6,\\
&\pi'_{19}=5\xi^5-3\xi^4+7\xi^3-24\xi^2-20\xi+6,\\
&\pi_{23}=-5\xi^5+4\xi^4-8\xi^3+26\xi^2+15\xi-9,\\
&\pi'_{23}=6\xi^5-4\xi^4+9\xi^3-30\xi^2-22\xi+6,\\
&\pi_{29}=11\xi^5-8\xi^4+17\xi^3-56\xi^2-35\xi+16,\\
&\pi_{31}=7\xi^5-5\xi^4+11\xi^3-36\xi^2-22\xi+7.
\end{align*}
$K=k(\sqrt{\eta})$ is a totally complex field of degree 12 over $\mathbb Q$ with the relative discriminant  $\mathcal{D}_{K/k}$ equal to $(\eta)$ as $\eta=\beta^2+4\gamma,$ where $\beta=\xi^5+\xi^4+\xi^3+1$, $\gamma=-173\xi^5+112\xi^4-270\xi^3+815\xi^2+576\xi-237$. From this we see that $d_K=7\cdot13\cdot19^2\cdot23^2\cdot29\cdot31\cdot23^2\cdot35509^2$. From Theorem~\ref{GT} we deduce that 
$$d(G_\varnothing)\geq t - r_1(k) - r_2(k) + \rho - 1=8-4-1+4-1=6.$$ The right hand side of the inequality from Theorem~\ref{inf} is equal to $2+2\sqrt{6}\approx6.8989<7,$ so it is enough to show that $d(G_T)>d(G_\varnothing),$ and to do this it is enough to construct a set of prime ideals $T$ and an extension of $K,$ ramified exactly at $T.$ 

Let $\pi_3=-6\xi^5+4\xi^4-9\xi^3+30\xi^2+21\xi-7$ be the generator of a prime ideal of norm 3 in $\mathcal{O}_k$ and $T$ be the set consisting of one prime ideal of $\mathcal{O}_K$ over $\pi_3\mathcal{O}_k.$ We see that $\pi_3\pi_{19}=11\xi^5-8\xi^4+17\xi^3-56\xi^2-35\xi+14=\rho^2+4\sigma,$ where $\rho=\xi^5+\xi^3+\xi^2+1,\, \sigma=2\xi^5-8\xi^4-14\xi^3-28\xi^2-9\xi+5,$ so $k(\sqrt{\pi_3\pi_{19}})/k$ is ramified exactly at $\pi_3$ and $\pi_{19}.$ But $\pi_{19}$ already ramifies in $K$ that is why $K(\sqrt{\pi_3\pi_{19}})/K$ is ramified exactly at $T$. Thus we have showed that $d(G_T)\geq 7$ and $K_T/K$ is indeed infinite.

To complete our proof we need a few more results.
\begin{theorem}[GRH Basic Inequality, see~\cite{Tsfa}, Theorem 3.1]
For an asymptotitically exact family of number fields under GRH one has:
\begin{equation}
\sum_q \frac{\phi_q \log q}{\sqrt{q}-1}+\phi_{\mathbb{R}}\left(\log 2\sqrt{2\pi}+\frac{\pi}{4}+\frac{\gamma}{2}\right)+\phi_{\mathbb{C}}(\log 8\pi+\gamma)\leq1,
\label{eqmain}
\end{equation}
the sum beeing taken over all prime powers $q.$
\end{theorem}
\begin{theorem}[See \cite{Haj1}, Theorem 1]
\label{estim}
Let $K$ be a number field of degree $n$ over $\mathbb Q,$ such that $K_T$ is infinite and assume that $K_T=\bigcup_{i=1}^{\infty} K_i.$ Then 
$$\lim_{i\to\infty} \frac{g_i}{n_i}\leq \frac{g_K}{n_K}+\frac{\sum_{\mathfrak{p}\in T}\log(N_{K/\mathbb{Q}}\mathfrak{p})}{2n_K}.$$
\end{theorem}

For our previously constructed field $K$ the genus is equal to $g_K\approx25.3490\hdots.$ From Theorem~\ref{estim} we easily see that $\phi_\mathbb{R}=0$ and $\frac{12}{2g_K+2\log 3}\leq\phi_\mathbb{C}\leq\frac{12}{2g_K},$ i.\:e., $0.23669<\phi_\mathbb{C}<0.22687.$ The lower bound for $\BS{K_T}$ is clearly equal to $$\BSL=1-\phi_\mathbb{R}\log2-\phi_\mathbb{C}\log(2\pi)\leq0.56498\hdots.$$ Knowing the decomposition in $K$ of small primes of $\mathbb{Q}$, we may now apply the linear programming approach to get the upper bound for $\BS{K_T}.$ This is done using the explicit formula (\ref{eqbs}) for the Brauer-Siegel ratio along with the basic inequality (\ref{eqmain}) and the inequality
$$\sum_{m=1}^{\infty} m\phi_{p^m}\leq \phi_{\mathbb R}+2\phi_{\mathbb C},$$
taken as the restrictions. This was done using the PARI package. As the calculations are rather cumbersome we will give here only the final result: $\BSU\approx0.59748\hdots,$ and the bound is attained for $\phi_7=\phi_9=\phi_{13}=0.03944\hdots, \phi_{19}=0.01002\hdots. \qed$
\section*{Acknowledgements}
I would like to thank my teacher professor M. A. Tsfasman, who gave the idea to study topics related to the Brauer-Siegel theorem and whose advices, suggestions and remarks were extremely valuable. I would also like to thank professor S. R. Louboutin for letting me know his results.


\begin{thebibliography}{20}
\bibitem{Bra}
R. Brauer `On zeta-functions of algebraic number fields', Amer. J. Math. 69, Num. 2, 1947, pp.243--250.
\bibitem{Gold}
D. M. Goldfeld `A simple proof of Siegel's theorem', Proc. Nat. Acad. Sci. USA 71 (1974), p. 1055.
\bibitem{Haj1}
F. Hajir, C. Maire `Tamely Ramified Towers and Discriminant Bounds for Number Fields II', Preprint.
\bibitem{Lan}
S. Lang `Algebraic number theory (Second Edition)', Graduate Texts in Mathematics 110, Springer-Verlag,
New York, 1994.
\bibitem{Lou00}
S. R. Louboutin `Explicit upper bounds for residues of Dedekind zeta functions and values of $L$-functions at $s=1$, and explicit lower bounds for relative class number of CM-fields', Canad. J. Math, Vol. 53(6), 2001, pp. 1194-1222.
\bibitem{Lou01}
S. R. Louboutin `Explicit lower bounds for residues at $s=1$ of Dedekind zeta functions and relative class numbers of CM-fields', Trans. Amer. Math. Soc. 355 (2003), 3079--3098.
\bibitem{Lou02} 
S. R. Louboutin. `On the Brauer-Siegel theorem.', J. London Math. Soc., to appear.
\bibitem{Mar}
J. Martinet `Tours de corps de classes et estimations de discriminants', Invent. Math. 44(1978), no. 1, 65-73.
\bibitem{Shaf}
I. Shafarevich `Extensions with prescribed ramification points', Publ. Math. I.H.E.S. 18(1964), pp. 71-95.
New York, 1994.
\bibitem{Sta}
H. M. Stark `Some effective cases of the Brauer-Siegel Theorem', Invent. Math. 23(1974), pp. 135--152.
\bibitem{Tsfa}
M. A. Tsfasman, S. G. Vl\u{a}du\c{t} `Asymptotic properties of global fields and generalized Brauer-Siegel Theorem', Moscow Mathematical Journal, Vol. 2, Num. 2, pp. 329-402.
\end{thebibliography}
\end{document}